\input amstex
\documentstyle{amsppt}
\input epsf.tex


\catcode`\@=11
\def\logo@{}
\catcode`\@=\active

\pageheight{46pc}

\mag\magstep 1
\NoBlackBoxes\tolerance10000

\def\ss{\smallskip}

\def\<{\langle}
\def\>{\rangle}

\def\inv{^{-1}}

\TagsAsMath
\def\ni{\noindent}
\def\Z{\Bbb Z}


\def\e{\epsilon}

\def\C{\Bbb C}




\def\ie{{{\it i.e.}\/}}

\def\({\hbox{\rm(}}
\def\){\hbox{\rm)}}

\tolerance 10000\NoBlackBoxes
\hoffset.2in

\mag\magstep 1
\def\ss{\smallskip}

\def\<{\langle}\def\>{\rangle}

\def\inv{^{-1}}

\TagsAsMath
\def\ni{\noindent}
\def\Z{\Bbb Z}

\def\C{\Bbb C}

\pageheight{46pc}
\leftheadtext{Gersten}
\rightheadtext{Asphericity}
\hoffset.2in

\mag\magstep 1
\def\ss{\smallskip}
\def\<{\langle}\def\>{\rangle}

\TagsAsMath
\def\ni{\noindent}
\def\Z{\Bbb Z}

\def\C{\Bbb C}

\topmatter
\title {Asphericity for certain groups of cohomological dimension~2}
\endtitle
\author{Steve Gersten}
\endauthor
\address{Mathematics Department, University of Utah, Salt Lake City, UT 84112}
\endaddress
\keywords{asphericity, cohomological dimension, idempotent,
Whitehead conjecture, Eilenberg-Ganea conjecture}
\endkeywords
\subjclass{20F06}
\endsubjclass
\email{sg\@math.utah.edu}
\endemail
\date{version 12, 3/16/2015}
\enddate
\thanks{Work supported by Social Security and Medicare,
\copyright{S. Gersten 2014}}
\endthanks
\thanks
{I am grateful to Jim Howie for his comments and suggestions.}
\endthanks

\abstract{A finite connected 2-complex~$K$ whose fundamental group~$G$
is of cohomological dimension~2 is aspherical 
iff the subgroup $\Sigma_K\subset H_2(K)$ of spherical
2-cycles is zero.

A finite connected subcomplex of an aspherical 2-complex is aspherical
iff its fundamental group is of cohomological dimension~2.
} 
\endabstract

\TagsAsMath
\address{Mathematics Department,
University of Utah, 
155 S. 1400 E., Salt Lake City, UT 84112-0090,
{\rm http://www.math.utah.edu/$\sim\null $sg}
}
\endaddress
\email{sg at math dot utah dot edu}
\endemail
\subjclass{20F06, 18F25}\endsubjclass

\endtopmatter
\document

\heading{\S 0. Introduction}\endheading

Two famous unsolved conjectures of low dimensional topology are the
Eilenberg-Ganea conjecture [EG] and the Whitehead conjecture [Wh].  The 
Eilenberg-Ganea conjecture states that a group of cohomological dimension  (cd)~2
is of geometric dimension (geomdim)~2 (the standard reference for these 
notions is K. Brown's book [Br] Chapter VIII \S 2).  The Whitehead conjecture
states that a connected subcomplex of an aspherical 2-complex is aspherical.
It is equivalent to saying that a connected subcomplex of a contractible
2-complex is aspherical.  There is an extensive
 literature on it including, notably,
work of J. F. Adams [Ad] and James Howie [Ho1].
A result of Bestvina and Brady [BB], now 17 years old, states that 
not both of these conjectures can be true. 
 If, for example, the Eilenberg-Ganea
conjecture is true, then the Whitehead conjecture must be false.  The proof
gives no indication about the truth or falsity of the individual conjectures.
There is also an article of Howie's [Ho2] which gives a new argument for
this result.

I view the Whitehead conjecture [Wh] as
 the combination of two distinct conjectures.
First, that a connected subcomplex~$K$ of an aspherical 
2-complex has fundamental
group~$G$ of cohomological dimension~2, and, second, the 
Eilenberg-Ganea conjecture,
that a group of cohomological dimension~2 has geometric dimension~2, so is
the fundamental group of an aspherical 2-complex.  This article deals with
the second conjecture in the case $K$ is finite.  I have nothing to say about
the first of these conjectures, since my methods apply to geometry only
when the fundamental group is already assumed to be of cohomological
dimension~2.

In this note I shall establish the following result.
\footnote{A stronger result is proved in \S 3 below after I
became aware of results of B. Eckmann [Ec].  I have 
not made any changes in the first two sections, but merely
added on the improvements at the end in \S 3.}

\proclaim{Theorem} Let $K$ be a finite connected 2-complex whose fundamental
group~$G$ is of cohomolgical dimension~2.
 Then $K$ is aspherical iff 
 \roster
 \item the relation module~$Z_1(\tilde K)$ is stably free, and
 \item the subgroup $\Sigma_K\subset H_2(K)$ of spherical
 2-cycles is 0.
 \endroster
 \endproclaim

The connection with the Whitehead conjecture is that every subcomplex~$K$ of
a contractible 2-complex has $H_2(K)=0$ (and hence, {\it a fortiori\/},
$\Sigma_K=0$), as one sees from the exact homology
sequence for the pair.  So, with the additional assumptions that the 
subcomplex is finite and connected and its fundamental group is of cd~2,
and the {\it necessary condition\/} (from the Theorem) 
that the relation module be
stably free,
it follows that $K$ is aspherical.
\demo{Remark}({\it due to Jim Howie}:)
Suppose $K$ is a  connected subcomplex of an {\it aspherical\/} 2-complex~$L$
(so we are replacing contractible by the wider context of aspherical).
Let $\bar K$ be the covering of $K$ contained in the universal cover~$\tilde L$
of $L$, $\bar K\subset \tilde L$, and let $\tilde K$ be the universal
cover of $K$.  So we have coverings $\tilde K\to \bar K\to K$, all with
the same $\pi_2$.  So the map $\pi_2( K)\to H_2(K)$ factors through
$H_2(\bar K)=0$.  It follows that $\Sigma_K=0$ in this situation,
whether or not $H_2(K)=0$.  
So with the further assumptions that $K$ is finite, that $\pi_1(K)$
is of cohomological dimension~2, and that the relation module
$Z_1(\tilde K)$ is stably free, it follows that $K$ is aspherical.
\enddemo

\heading{\S 1.  Proof of Theorem and some Corollaries
\footnote{See \S 3 for the strengthened version of the Theorem of \S 0.}
}
\endheading

In this section $K$ is a finite connected 2-complex whose fundamental group~$G$
is of cohomological dimension~2,
and $\tilde K$ is the universal covering complex of $K$.  We have 
the fibration $G\to \tilde K\to K$ with discrete fibre $G$.  This is
classified\footnote{See any good book on algebraic topology,
like [Sp].} by a map $K\to BG$ with fibre $\tilde K$, where $BG$ is an
Eilenberg-MacLane space of type $K(G,1)$.  We can apply the Serre
spectral sequence to this fibration to get 
the $E^2_{p,q}$ term $H_p(G,H_q(\tilde K),\Z)$ converging to the $E^\infty$
term, which is the successive quotients of a filtration of $H_{p+q}(K,\Z)$.

The the spectral sequence
has only two non-zero rows, rows~0 and~2, since $\tilde K$ is simply
connected (so row~1 is zero) and since $\tilde K$ is a 2-complex (so
rows~$n$ with $n\ge 3$ are zero).
So we have an edge-term exact sequence
$$0\to H_3(G,\Z)\to H_0(G,H_2(\tilde K,\Z))\to H_2(K,\Z)\to H_2(G,Z)\to 0.$$

The 0 on the left comes about since $K$ is a 2-complex so $H_3(K)=0$.
Since cd($G$)=2, it follows that $H_3(G,\Z)=0$, so
the exact sequence simplifies to
$$0\to  H_0(G,H_2(\tilde K,\Z))\to H_2(K,\Z)\to H_2(G,Z)\to 0.\tag 1.1$$

\proclaim{Lemma 1.2}  If\ $\Sigma_K=0$, then $H_0(G,H_2(\tilde K,\Z))=0.$
\endproclaim

For the proof, note that the 
image of the map $H_0(G,H_2(\tilde K,\Z))\to H_2(K,\Z)$
is $\Sigma_K$, the group of spherical 2-cycles of $K$.

Now we examine the hypothesis that cd($G$)=2 more closely.
Let $P=: H_2(\tilde K,\Z)$.

\proclaim{Lemma 1.3}  
$P$ is a finitely generated
projective $\Z [G]$-module and a direct summand of $C_2(\tilde K,\Z)$.
\endproclaim

\demo{Proof}
Consider the exact sequence
$$0\to Z_1(\tilde K)\to C_1(\tilde K)\to C_0(\tilde K)\to\Z\to 0.\tag 1.4$$
But $C_1$ and $C_0$ are free $\Z [G]$-modules and hence projective.
So it follows from [Br] Chapter VIII \S 2 Lemma 2.1 (iv) (using $n=2$
there and the fact that cd($G$)=2)
that $Z_1(\tilde K)$ is projective.

Now $B_1(\tilde K)=Z_1(\tilde K)$ since $\tilde K$ is simply connected.
It follows that the
exact sequence
$$0\to H_2(\tilde K)\to C_2(\tilde K)\to B_1(\tilde K)\to 0\tag 1.5$$
splits.
Thus $P$ is projective and a direct summand of $C_2(\tilde K)$ with
complementary direct summand $Z_1(\tilde K)=B_1(\tilde K)$,
and the lemma is proved.
\enddemo

Now assume that $K$ is aspherical.  Then
$P=H_2(\tilde K)=0$, so $C_2(\tilde K)\cong B_1(\tilde K)=Z_1(\tilde K)$,
whence $Z_1(\tilde K)$ is free (and hence {\it a fortiori\/}
stably free).

Summarizing, if $K$ is aspherical, it follows that $\Sigma_K=0$
and $Z_1(\tilde K)$ is stably free, which proves the necessity of
conditions (1) and (2) in the Theorem.

Next we proceed to the sufficiency of conditions (1) and (2).
Since $\Sigma_K=0$, it follows from Lemma~1.2
that 
$H_0(G,H_2(\tilde K,\Z))=0.$  But this is the same as $P/IP$,
where $P=H_2(\tilde K,\Z))$ and $I$ is the augmentation ideal
of $\Z[G]$.
Thus $P/IP=0$.

Furthermore we are assuming that $Z_1(\tilde K)$ is stably free.
We do elementary expansions on $K$, attaching a finite number of
2-discs trivially at the base point.  The effect on the chains
of the universal cover is to take the direct sum of $C_2$ and $C_1$ with
the same finite number of copies of $\Z[G]$, and map them by the identity
map.  So $Z_1=B_1$ is increased by the direct sum of the same number
of copies of $\Z[G]$.  This has the effect of stabilizing $Z_1$, which
we may assume now is free.  In addition there is no change in 
$H_2$ of $\tilde K$, which is still $P$.
So we may, without loss of generality, adjust the notation and
assume that $K$ is such that $Z_1(\tilde K)$ is free and
$P/IP=0$, with $P=H_2(\tilde K)$.
Our goal is to prove that $P=0$.

Now take the split exact sequence 1.4 and apply the functor $\Z\otimes_{\Z[G]}$
to it.
The left term is $P/IP$ and hence vanishes, so
the result is
that $\Z\otimes_{\Z[G]}C_2(\tilde K)\cong \Z\otimes_{\Z[G]} Z_1(\tilde K)$.
But both $C_2$ and $Z_1$ are free $\Z[G]$ modules, so it follows that
they are of the same rank.
Thus the boundary map $\partial:C_2(\tilde K)\to Z_1(\tilde K)$
is a split epimorphism of free modules of the same rank.
Let $h: Z_1(\tilde K)\to C_2(\tilde K)$ be a right inverse for $\partial$,
so $\partial\circ h=Id$, the identity map of $Z_1$.  Since
$C_2$ and $Z_	1$ are free of the same rank and equipped
with free bases, the maps
$\partial$ and $h$ are given by $n$-by-$n$ matrices
over $\Z[G]$ with respect to these bases, so we are
dealing with $M_n(\Z[G])$.

Let $E=:h\circ \partial:C_2(\tilde K)\to C_2(\tilde K)$.
Computation shows that $E^2=E$, so the image of $E$ is a direct summand
of $C_2(\tilde K)$ which is isomorphic to $Z_1(\tilde K)$.  It follows
that $Id-E$ projects $C_2(\tilde K)$ onto an isomorphic copy of $P$,
where here $Id$ is the identity map of $C_2$. 

Let $\Lambda$ be the ring of $\Z[G]$-linear endomorphisms of $C_2(\tilde K)$,
so $\Lambda\cong M_n(\Z[G])$ where $n$ is the rank of $C_2(\tilde K)$ as a free
$\Z[G]$-module.  We consider $M_n(\Z[G])=M_n(\Z)[G]$, the group ring
of $G$ with coefficients in $M_n(\Z)$, and extend scalars to $\C$ to
get $M_n(\C)[G]$.  A typical element of this ring
may be written as $z=:\sum_g m_g g$, where $m_g\in M_n(\C)$
and almost all $m_g$ are zero.

 Define \footnote{I am following almost {\it
 verbatim\/} the exposition given in [Mo], with appropriate changes for the
 wider context of $n$-by-$n$ matrices rather than 1-by-1.}
  a map $t:M_n(\C)[G]\to \C$
by $t(z)=: {\text{trace}}(m_1)$, where $m_1$ is the coefficient of 1, the unit
element of $G$, in $z$.  
In addition $M_n(\C)[G]$ has an involution $x\to x^*$ given by
$x^*=\sum_{g\in G}\bar {x_g}^T g\inv$ if $x=\sum_{g\in G}  x_g g$;
here the superscript~$T$ indicates the transposed matrix and 
``bar" is complex
conjugation.
The following facts are valid for $t$.
\roster
\item $t(Id)=n$, where $Id$ is the identity matrix.
\item $t$ is $\C$-linear.
\item $t(xy)=t(yx)$ for all $x,y$ in $M_n(\C)[G]$.
\item $t(xx^*)\ge 0$ and $t(xx^*)=0$ iff $x=0$.
\endroster
Only the last item requires comment.
Let $x=\sum_{g\in G}x_g g$ and let $\bold r_i(x_g)$ denote the
$i^{\text{th}}$ row of the matrix $x_g$.
Then $t(xx^*)=\sum_{g\in G}\sum_{i=1}^n||\bold r_i(x_g)||^2$,
where, for a row vector $\bold v$, $||\bold v||$ is the $L^2$-norm of the
vector.\footnote{So here $\bold r_i(x_g)$ is the vector over $\C$
 with components
$(x_{g,i1},x_{g,i2},\dots,x_{g,in})$ and $||\bold r_i(x_g)||^2=
\sum_k|x_{g,ik}|^2$.
}
Now $M_n(\C)[G]$ is equipped with a Hermitian inner product given
by $<\sum x_g,\sum y_g>=t(\sum x_g {\bar y_g}^T)$.  The non-degeneracy
follows from the identity $<x,x>=t(xx^*) \ge 0$ and equal to 0
iff $x=0$.  It follows that $M_n(\C)[G]$ may be completed to a
Hilbert space~$H$.  Elements of $\Lambda$
 act as bounded operators $B(H)$ on $H$.
Let $A$ denote the closure of the set of these operators in the norm
topology of $B(H)$.  Then $A$ is a $C^*$ algebra, so $1+ff^*$ is invertible
for all $f\in A$.
\footnote{See W. Rudin, Functional Analysis, Second Edition,
McGraw-Hill, Inc., 1991, p. 295, Theorem 11.28 (f).
}
Recall that an element~$x$ is self-adjoint if $x=x^*$ and a projection
if it is  both idempotent and self-adjoint.  The next two results are
due to Kaplansky.  I am following [Ma] in the exposition.

\ss

\ni 1.5 Let $R$ be a unital ring with involution~${\null}^*$ such that for any
$x\in R$ the element $1+xx^*$ is invertible in $R$.  For every
idempotent $f\in R$ we can find a projection $e\in R$ so that
$Rf=Re$.  (For the proof see [Ma] Prop.~2.4.  Note that we are using
left modules rather the right modules in [Ma].)

\ss

\ni 1.6 Let R be a unital ring and $e,f\in R$ two idempotents such 
that $Re=Rf$.  Then $e$ and $f$ are similar, so
there exists an invertible element~$s\in R$ so that
$s\inv f s=e$.  (See [Ma] Prop.~2.5.)
\ss
It follows that every idempotent in $A$ is similar to a projection.
So the idempotent $E\in M_n(\Z[G])$ is similar in the larger ring
to a projection $\Pi$.
Also the function $t$ is continuous so extends a function also called
$t$ on A.
I also recall Lemma~1 of [Mo] which asserts that if $t(xx^*)=0$ for
$x\in A$, then $x=0$.  Although the context is more restricted
there, the same argument
given there works here (it is, basically, a continuity argument).

\proclaim{Lemma 1.7} For the idempotent $E$ in $M_n(\C)[G]$,
$t(E)$ is real and $0\le t(E)\le n$.
\endproclaim

\demo{Proof}  $E$ is similar in $A$ to the projection $\Pi$,
so $t(E)=t(\Pi)$, using property 3.  But $\Pi=\Pi\cdot \Pi^*$,
so $t(\Pi)>0$ if $\Pi\ne 0$, by the analog of Lemma~1 of [Mo].
Similarly $t(Id-\Pi)>0$ if $\Pi\ne Id$.
\enddemo 

Now we can complete the proof of the Theorem.  We have $\partial\cdot h$
is the identity of $Z_1(\tilde K)$, a free module of rank~$n$,
so $t(\partial\cdot h)=n$.
But $t(\partial \cdot h)=t(h\cdot \partial)$, so $t(E)=t(h\cdot \partial)=n$.
Thus $t(Id-E)=0$, and it follows that $Id-E=0$.  Hence $E$ is the identity,
the maps $\partial$ and $h$ are mutual inverses, $P=0$,
$K$ is aspherical, and the proof is complete.\qed

Let us deduce some corollaries of the Theorem.
\proclaim{Corollary 1}  If $K$ is a finite connected 2-complex whose
fundamental group is of cohomological dimension~2 and if $H_2(K)=0$ 
and the relation module~$Z_1(\tilde K)$ is stably free,
then $K$ is aspherical.  In particular, if $K$ is a finite connected
subcomplex of a contractible 2-complex, and if the cohomological 
dimension of $\pi_1(K)$ is~2
and if the relation module is stably free, then $K$ is aspherical.\qed
\endproclaim

\proclaim{Corollary 2}  Suppose $G$ is a finitely presented perfect group
\footnote{Or more generally if its first betti number is 0.}
of cohomological dimension~2 which admits a balanced
\footnote{A finite presentation is balanced if it has the same number
of generators as defining relators.} presentation~$\Cal P$.
  If $K=K(\Cal P)$ is the associated 
2-complex to the presentation
\footnote{$K(\Cal P)$ has a single vertex, one 1-cell for
each generator, and one 2-cell for each defining relator
whose attaching map spells out the relator in the 1-skeleton,
with appropriate orientations.}
 and if the relation module is stably free, then $K$ is aspherical.
\endproclaim

\demo{Proof}  An Euler characteristic computation shows that the
second betti number of $K$ is 0.  Since $K$ is a 2-complex, this implies
$H_2(K)=0$, so Corollary~1 applies.
\enddemo

\demo{Remark}  Perfect groups are in a sense opposite to locally
indicable groups.  Both Howie's and Adams's methods [Ho1],[Ad] fail
for perfect groups, and Howie's apply to locally indicable groups.
So Corollary~2 above reaches into territory not explored in [Ho1]
at the cost of having to assume cohomological dimension~2 and the ubiquitous
necessary condition that the relation module be stably free.
\enddemo

\proclaim{Corollary 3}  Let $K$ be a finite connected 2-complex whose
fundamental group $G$ is of cohomological dimension 2
and such that $G$ is of type FL.
\footnote{A group $G$ is of type FL if the trivial $G$-module
$\Z$ has a finite resolution by finitely generated free
$\Z[G]$-modules.
} 
Suppose in addition that
$H_2(K)=0$.  Then K is aspherical.
\endproclaim

\demo{Proof}  One checks that conditions 1 and 2 in
the Theorem are satisfied, with 2 following since there
are no nonzero 2-cycles.
Since $G$ is of type FL, there is a finite resolution
$$0\to L_n\to L_{n-1}\to ... \to L_1\to L_0\to \Z,$$
with the $L_i$ finitely generated free $\Z[G]$-modules.
Let $Z_i$ be the kernel of the map $L_i\to L_{i-1}$.
Since $G$ is assumed of cohomological dimension~2,
$Z_1$ is finitely generated and projective.  It follows
that all the short exact sequences
$$0\to Z_{i}\to L_{i}\to Z_{i-1}\to 0$$
split for $2\le i\le n$, where $Z_n=L_n$.  
Now it follows by downward induction that 
$Z_n, Z_{n-1},\dots, Z_2,Z_1$ are stably free.
In particular $Z_1$ is stably free.  But $Z_1$ is
stably equivalent to the relation module for $K$,
namely, $Z_1(\tilde K)$, and it follows that 
$Z_1(\tilde K)$ is stably free.

Hence both conditions of the Theorem are satisfied,
and it follows that $K$ is aspherical.
\enddemo

There is a version of the last corollary which is analgous
to Howie's observation in \S 0, as follows.

\proclaim{Corollary 4}  Let $K$ be a finite connected
subcomplex of an aspherical 2-complex~$L$.  Assume
that the fundamental group~$G$ of $K$ is of cohomological
dimension~2 and of type FL.  Then $K$ is aspherical.\qed
\endproclaim

\demo{Remark}  Corollary 4 is in fact a necessary and sufficient 
condition for a finite connected subcomplex $K$ of an aspherical 2-complex $L$
to be aspherical, namely, that $G=\pi_1(K)$ be of cohomological
dimension (at most) 2 and that it be of type FL.
The sufficiency is the content of Corollary~4.  For necessity,
assume that such $K$ is aspherical.  Then the chain complex of 
$\tilde K$ gives a finite free resolution of $\Z$ by finitely
generated free $\Z[G]$-modules, so $G$ is of type FL.  The fact that
$G$ is of cohomological dimension (at most) 2 follows from
[Br] p. 184, Lemma~2.1(iv).
\enddemo

\heading{\S 2. The {\rm t}-rank of a projective module
and some open questions}
\endheading

\demo{Definition}  The t-rank of $P$, a finitely generated projective
module over $\Z[G]$, is defined to be $t(E)$, where $E$ is an idempotent
matrix in $M_n(\Z)[G]$ with image~$P$.  For simplicity denote
the t-rank of $P$ by $t(P)$.
\enddemo

Here is a collection of properties of the t-rank, all easily proved.
\roster
\item $t(P)$ is well-defined and independent of the ambient matrix algebra
used to define it.
\item $t(P\oplus Q)=t(P)+t(Q)$.
\item If $F$ is a free module over $\Z[G]$ of rank~$n<\infty$, then
$t(F)=n$.\footnote{This makes sense since $\Z[G]$ has invariant
basis number.}
\item $t(P)$ is a non-negative integer with 
 $t(P)=0$ iff $P=0$.
\endroster

There is another notion of rank of a f.g. projective module~$P$
over $\Z[G]$, namely the rank over $\Z$ of
$Z\otimes_{\Z[G]}P$.  
It is easy to see that the two notions of rank agree on free
modules, and more generally they agree on stably free 
modules.

\demo{Remark}  Let $E$ be an idempotent matrix in
$M_n(\Z[G])$ with image the projective module~$P$.
Then the two notions of rank are easily computed as
follows.  If $E=\sum_{g\in G}   e_{g}g$ where $e_g\in M_n(\Z)$,
and if $e_g$ is the matrix $(e_{g,ij})$, then the $t$-rank
is given by
$t(P)=\sum_{1\le i\le n} e_{1,ii}$ whereas 
the other notion of rank is calculated as
\footnote{Over $\Z$, the target of $\e$, any idempotent
matrix is similar to a diagonal matrix with 1's and 0's on
the main diagonal, and the rank is the number of 1's,
\ie\ the trace.  Over
the commutative ring $\Z$ trace is a similarity invariant,
from which the rank formula can be derived.}
$\e(\text{trace}(E))=\sum_{g,i}e_{g,ii}$.  If $P$ is stably free,
the two calculations must give the same answer.
\enddemo

\demo{Question}  Do the two notions of rank agree on all
finitely generated projective modules over $\Z[G]$?
If there is an projective module~$P$ for which the two notions 
disagree, that
module cannot be stably free.
In particular, if $P/IP=0$ but $P\ne 0$, then $P$ cannot
be stably free.
\enddemo

\demo{Question}  If $P$ is a finitely generated projective module over an integral group ring such that $P/IP=0$, is $P=0$?  
\footnote{
I was motivated originally to ask this
question by Nakayama's lemma, but it
has turned out to be much deeper.  In a sense it says that
the augmentation ideal plays the role of the Jacobson radical
as far as projective modules over group rings are concerned.
I had originally approached the question, whether $P/IP=0$ implies $P=0$,
as a question of ring theory.  Howie pointed out to me that there is a
canonical isomorphism $I/I^2\cong G_{ab}$, where $G_{ab}$ is the abelianization
of the group~$G$ and I is the augmentation ideal of $\Z[G]$.  Hence if $G$
is perfect, $I=I^2=I^3=\dots =I^\infty$, and the ring theoretic approach
fails.  But for abelian groups and potentially a larger class admitting
no perfect subgroups, this approach may conceivably succeed.
}
\enddemo

\demo{Question}  Is every relation module of a finite
presentation of a group~$G$ of cohomological dimension~2
stably free?  It is {\it{equivalent\/}} to ask whether such
$G$ is of type FL.
\enddemo

\demo{Question}  If $P$ is a finitely generated projective
module over $\Z[G]$, is $\ell_2(G)\otimes_{\Z[G]}P$ stably free
\footnote{Putatively it should  be free of rank~$t(P)$.}?
An affirmative answer to this question would imply affirmative
answers to the preceding two questions.
\footnote{I have asked the experts in the field of operator
algebras and come up with either no response or that they don't
know the answer.}
\enddemo

\demo{Remark}
I want to pose a vague question to bring the Eilenberg-Ganea
conjecture into the discussion: what in addition to 
cohomological dimension~2 and type FL 
(both necessary conditions) does one need to 
characterize when a finitely presented group has a
finite aspherical
presentation?  
\footnote{I have tried and failed so far to prove that just these
two conditions sufficed.  I am grateful to Howie for having 
found the gap in my argument.}
\enddemo

\demo{Remark}  My results have {\it no bearing\/} on the result
of [BB], that one of the two conjectures, Whitehead's
or Eilenberg-Ganea's, must be false.  The reason is I am
considering only finite presentations.  The Whitehead conjecture
considers arbitrary {\it not necessarily finite\/} subcomplexes
of a contractible 2-complex and the Eilenberg-Ganea
conjecture considers arbitrary groups, {\it not necessarily
finitely presented\/}, of cohomological dimension~2.  In this
generality, one of the two conjectures must be false.
\enddemo

\heading{\S 3. Update and further results}
\endheading

Since sections 0--2 were written, I have become aware of an article
by B. Eckmann [Ec] whose results enable me to strengthen the Theorem
of \S 1 and to settle some of the questions asked at the end of \S 2,
at least for groups of cohomological dimension~2.
Eckmann's results for our purposes can be summarized as follows.
\footnote{The statement that follows is a combination of several of his
results, where I have omitted references to various forms of the 
Bass conjecture to limit demands on the reader.  Also I only need the results for groups of cd 2 over
$\Z$, whereas he considers other fields.}

\proclaim{Theorem of Eckmann}
\footnote{I am thankful to Stefan Friedl for catching an omission
in the statement
in an earlier version of this paper.
}  If $G$ is a countable group of cohomological
dimension~2 over $\Z$ and $P$ is any finitely generated 
projective $\Z[G]$-module,
then $\ell_2(G)\otimes_{\Z[G]}P$ is a free $\ell_2(G)$-module of
rank the $\Z$-rank of $\Z\otimes_{\Z[G]}P$.
\endproclaim

To apply this result we need a lemma, which makes use of
 the
central point of the proof of the Theorem of \S 1.

\proclaim{Lemma 3.1} Let $G$ be a group with the property that
for any finitely generated projective module $P$ over $\Z[G]$
one has $\ell_2(G)\otimes_{\Z[G]}P$ is a free $\ell_2(G)$-module of
rank the $\Z$-rank of $\Z\otimes_{\Z[G]}P$.
If $Q$ is a finitely generated projective $\Z[G]$-module such that
$Q/IQ=0$, then $Q=0$.
\endproclaim

\demo{Proof}  Let $P$ be a finitely generated projective such that
$Q\oplus P=\Z[G]^n$.
Since $Q/IQ=0$, it follows that the $\Z$-rank of $\Z\otimes_{\Z[G]}P=n$.
Now tensor with $\ell_2(G)$ over $\Z[G]$ to get $\ell_2(G)\otimes_{\Z[G]}P$
is free of rank~$n$ as an $\ell_2(G)$-module.  
It follows that we have homomorphisms of free $\ell_2(G)$-modules
 $\partial :\ell_2(G)^n\to \ell_2(G)\otimes_{\Z[G]}P$
 and $h:\ell_2(G)\otimes_{\Z[G]}P\to\ell_2(G)^n$ of the same rank~$n$
 such that
 $\partial\circ h=1_{\ell_2(G)\otimes_{\Z[G]}P}$.
 Since these maps are given by $n$-by-$n$ matrices over $\ell_2(G)$,
 we can apply the function $t$ (which extends by continuity to
 $M_n(\ell_2(G))$) to get
 $n=t(1_P)=t(\partial \circ h)=t(h\circ \partial)$.  But
 $h\circ \partial$ is an idempotent endomorphism of $\ell_2(G)^n$.
 It follows that $t(\ell_2(G)\otimes_{\Z[G]}Q)=t(Id - h\circ \partial)=n-n=0$,
 where $Id$ is the identity map of $\ell_2(G)^n$.
 Hence $\ell_2(G)\otimes_{\Z[G]}Q=0$, whence $Q=0$, and the proof is 
 complete.
\enddemo

\proclaim{Corollary 3.2}  If $G$ is a countable group of cohomological
dimension~2 over $\Z$ and $P$ is a
 finitely generated projective module over $\Z[G]$
such that $P/IP=0$, then $P=0$.
\endproclaim

\demo{Proof} 
This is immediate from Eckmann's theorem and the lemma.
\enddemo

We can now strenthen the Theorem of \S 1.

\proclaim{Main Theorem 3.3}  Let $K$ be a finite connected 2-complex whose
fundamental group is of cohomological dimension~2.
Assume that the subgroup~$\Sigma_K\subset H_2(K)$ of spherical 2-cycles
is 0.  
Then it follows that
\roster
\item $K$ is aspherical, and
\item the relation module $Z_1(\tilde K)$ is stably free, and
\item the the fundamental group $G$ of $K$ is of type FL.
\endroster
\endproclaim

\demo{Proof} From the exact sequence 1.1 we deduce from the
hypothesis that $\Sigma_K=0$ that 
$H_0(G,H_2(\tilde K,\Z)=0$.  That is, if we let $P=H_2(\tilde K,\Z)$,
$P/IP=0$.  Since $P$ is finitely generated and projective over $\Z[G]$,
it follows from Corollary~3.2 that $P=0$.
Thus $K$ is aspherical.
Also from the split exact sequence 
$$0\to H_2(\tilde K)\to C_2(\tilde K)\to Z_1(\tilde K)\to 0$$
we deduce that $Z_1(\tilde K)$ is stably free.
It follows easily that $G$ is of type FL.
This completes the proof.
\enddemo

\proclaim{Corollary 3.4}  Let $K$ be a finite connected subcomplex of an
aspherical 2-complex such that the fundamental group~$G$ of $K$ is
of cohomological dimension~2.  Then $K$ is aspherical.
\endproclaim

\demo{Proof}  This is immediate from Corollary~4 of \S 2 and 
Theorem~3.3.
\enddemo

\demo{Remark}  The question remains open whether $P/IP=0$ implies
$P=0$ for a finitely generated projective module over the group
ring~$\Z[G]$ of a general group~$G$.  There are many classes of
groups for which this holds by our argument and by the results
of the appendix to [Ec]; the question is intimately related to
the various forms of the Bass conjecture, which we have chosen not
to discuss here.
\enddemo

\heading{\S 4. Some decision problems for finitely presented
groups of cohomological dimension~2}
\endheading

The Main Theorem of \S 3
reduces the question of asphericity of $K=K(\Cal P)$ to
a question about $\Sigma_K\subset H_2(K)$, when $\Cal P$
is a finite presentation of a group $G$ with $\text{cd}(G)=2$.
Thus there is the question about how effective are calculations
involving $\Sigma_K\subset H_2(K)$ and $H_2(G)=H_2(K)/\Sigma_K$.  

For general finitely presented groups~$G$ there is no effective procedure
for calculating $H_2(G)$ nor even for finding its rank nor whether or not
it is zero [Go].  The situation is different for finitely presented groups
of cohomological dimension~2, but the statement of the result depends
on how the data are given.

Let $K=K(\Cal P)$ be a the 2-complex associated to a finite presentation
$\Cal P$ of the group~$G$ of cohomological dimension~2.  We have seen
(using the fact that $B_1(\tilde K)=Z_1(\tilde K)$)
that $Z_1(\tilde K)$ is a finitely generated projective $\Z[G]$ module,
and the short exact sequence
$$0\to H_2(\tilde K)\to C_2(\tilde K)\overset{\partial}\to
\rightarrow Z_1(\tilde K)\to 0$$
splits.  Thus there is a section $h:Z_1(\tilde K)\to C_2(\tilde K)$
so that $\partial\circ h=1$, the identity map of $Z_1(\tilde K)$.
It follows that $E=:h\circ\partial$ is an idempotent endomorphism 
of $C_2(\tilde K)$; the kernel of $E$ is isomorphic to
$H_2(\tilde K)$ and the image of $E$ is
$Z_1(\tilde K)$.  To avoid questions about the effectiveness of choosing
$h$, we shall assume that $E$ is given, say, as an n-by-n
matrix with entries in $\Z[G]$, where $n$ is the rank of $C_2(\tilde K)$
as a free $\Z[G]$-module.
\footnote{The various definitions of projective module are all equivalent
in classical mathematics, but some are more effective than others.
We have chosen to formulate the
proposition in terms of the most effective definition by idempotent
matrices.}

\proclaim{Proposition 4.1}  In the situation just described, assume
$E$ is given as an idempotent matrix over $\Z[G]$ whose image
is isomorphic to
$Z_1(\tilde K)$ and such that the kernel of $E$ is isomorphic to
$H_2(\tilde K)$.
Then the following are true.
\roster
\item $\Sigma_K$ is effectively computable,
\item $H_2(G)$ is effectively computable,
\item $\beta_2(G)$ is effectively computable, and
\item one can effectively decide whether or not $H_2(G)$ is zero.
\item $\Sigma_K$ is a direct summand of $H_2(K)$
and $H_2(G)$ is a torsion-free abelian group.
\endroster
\endproclaim

\demo{Remark}  It is convenient to factor $\partial: C_2(\tilde K)\to C_1(\tilde
K)$ as $\partial =\iota\circ p$, where $p:C_2(\tilde K)\to B_1(\tilde
K)=Z_1(\tilde K)$ is a split epimorphism
of finitely generated projective $\Z[G]$-moudles
 and $\iota:Z_1(\tilde K)\to C_1(\tilde K)$
is injective.  Then $H_2(K)=\text{ker}(C_2(K)\to C_1(K))$ and
$\Sigma_K=\Z\otimes_{\Z[G]} \text{ker}(E)$.  Since $E=h\circ p$,
one has $(\iota\circ p)\circ E=\iota\circ p\circ ( h\circ p)=\iota\circ p$,
and one sees that $\Sigma_K\subset H_2(K)$, as a check on the notation.
\enddemo

\demo{Proof of 4.1}  One has $\Sigma_K=\text{ker}(\e(E))$, where
we recall that $\e:\Z[G]\to \Z$ is the augmentation.  
Similarly $H_2(K)$ is the kernel of a square matrix over $\Z$.
It follows that conclusions 1--4 are consequences of the theorem
of elementary divisors for matrices over $\Z$.

For the final conclusion I show first that $H_2(G)$ is torsion-free.
\footnote{I thank Jim Howie for this argument.}  
We have for any prime number~$p$,
 $0=H^3(G,\Z/p\Z) =\text{Hom}(H_3(G,\Z),\Z/p\Z)\oplus 
\text{Ext}(H_2(G,\Z),\Z/p\Z)=\text{Ext}(H_2(G,\Z),\Z/p\Z)$,
by the universal coefficient theorem and the assumption that $\text{cd}(G)=2$.
But if there were any non-zero torsion in $H_2(G)$, it would show up
in the Ext terms.  It follows that $H_2(G)$ is torsion-free.

Finally, since $H_2(G)=H_2(K)/\Sigma_K$ is torsion-free, it follows that
all of the elementary divisors for the pair $(H_2(K),\Sigma_K)$ are 
either 0 or 1.  Hence $\Sigma_K$ is a direct summand of 
$H_2(K)$, and the proof is complete.
\enddemo

\demo{Problem}  Can one by applying a
sequence Tietze transformations to $\Cal P$ arrange that $\Sigma_K=0$?
This is just the special case of the Eilenberg-Ganea conjecture for
finitely presented groups of cohomological dimension~2.
\footnote{Conclusion 5 of the proposition can be interpreted as saying
there are no obvious obstructions to this case of the Eilenberg-Ganea
conjecture being true.}
\enddemo

\end
\bye